\documentclass[a4paper]{amsart}
\usepackage{amsmath, amssymb}
\newtheorem{Theorem}{Theorem}[section]

\theoremstyle{definition}

\theoremstyle{remark}

\numberwithin{equation}{section}
\newcommand{\R}{\mathbb R}
\newcommand{\N}{\mathbb N}

\newcommand{\A}{\mathcal{A}}

\newcommand{\s}{\mathcal{S}}
\begin{document}

\title{From infinitesimal symmetries to deformed symmetries of Lax-type equations}
\author{Jean-Pierre Magnot}
%\author{}%

\address{Lyc\'ee Jeanne d'Arc - Avenue de Grand Bratagne - F-63000 Clermont-Ferrand}
\email{jean-pierr.magnot@ac-clermont.fr}
%\email{}
%\email{}%

%\commby{}%
% ----------------------------------------------------------------

\begin{abstract}
Using the procedure initiated in \cite{Ma2013}, we deform Lax-type equations though a scaling of the time parameter. This gives an equivalent (deformed) equation which is integrable in terms of power series of the scaling parameter. We then describe a regular Frölicher Lie group of symmetries of this deformed equation      
\end{abstract}
\maketitle


\section*{Introduction}
In \cite{Ma2013}, an algebra and a group of formal series of operators is described in order to rewrite the integration of the KP hierarchy in a non formal way. One of the main advances of this work is to get a (non formal) principal bundle where the concept of holonomy makes sense rigorously. The geometric objects under consideration are diffeological or Fr\"olicher groups, which are regular in the sense that the exponential map exists and is smooth. Diffeological spaces, first described in the 80's by Souriau and his coworkers \cite{Don, Igdiff, Les, Sou} are generalizations of manifolds that enables differential geometry without charts. Independently, Fr\"olicher spaces give a more rigid framework, that also generalize the notion of manifolds \cite{CN, FK, KM}. The comparison of the two frameworks has been made independently in \cite{Ma2006} and in \cite{Wa}, see e.g. \cite{Ma2013}. The aim of this paper is to show how this framework can apply to the theory of Lax equations. A Lax equation \cite{Lax} is a formally integrable equation of the type : 
$\partial_t L = [P,L]$
where $P,L$ are in most cases differential, pseudo-differential, or difference operators. This equation integrates heuristically  as a classical equation on a group of matrices: there should have an unique solution, up to the initial value $L(0)$, given by $L(t) = Ad_{Exp P (t)} L(0).$
Unfortunately, very often, the operator   $ Exp P (t)$ exists only at a formal level. We propose to apply a scaling $t \mapsto qt$ to the time variable. The operator $P(t)$ is changed into an operator $P_q(t)$ which is a monomial of order $1$ in the $q-$variable, adapting the ideas of \cite{Ma2013}. This allows the machinery of $q-$deformed operators: the algebras considered are now Lie algebras of (smooth) regular Lie groups. As a simple consequence, we get smoothness of the unique solution $L_q(t)$ with respect to $P(t)$ and $L(0);$ another consequence is that the full space of symmetries is a diffeological or Fr\"olicher group, and that a class of symmetries of the $q-$deformed Lax equation obey also a Lax-type equation $\partial_tS_q = [ad_{P_q}, S_q] .$ This equation is here interpreted as a holonomy equation, which integrates by virtue of the results of \cite{Ma2013}. These symmetries $S_q$ are rigorously constructed when $L_q$ is a q-deformed formal pseudo-differential operator.
\section{Preliminaries: Fr\"olicher Lie groups of formal series}

We now turn to  key results from \cite{Ma2013}: 

\begin{Theorem} \label{regulardeformation}
Let $(A_n)_{n \in \N^*} $ be a sequence of complete locally convex (Fr\"olicher)
vector spaces which are regular, 
equipped with a graded smooth multiplication operation
on $ \bigoplus_{n \in \N^*} A_n ,$ i.e. a multiplication such that 
$A_n .A_m \subset A_{n+m},$ smooth with respect to the corresponding Fr\"olicher structures.
Then, the set 
$1 + \A = \left\{ 1 + \sum_{n \in \N^*} a_n | \forall n \in \N^* , a_n \in A_n \right\} $
is a Fr\"olicher Lie group, with regular  Fr\"olicher Lie algebra
$\A= \left\{ \sum_{n \in \N^*} a_n | \forall n \in \N^* , a_n \in A_n \right\}.$
Moreover, the exponential map defines a bijection $\A \rightarrow 1+\A.$  
\end{Theorem}
%\begin{Theorem}\label{exactsequence}
%Let 
%$$ 1 \rightarrow K \underrightarrow{i} G \underrightarrow{p}  H \rightarrow 1 $$
%be an exact sequence of Fr\"olicher Lie groups, such that there is a smooth section $s : H \rightarrow G,$ and such that 
%the trace diffeology from $G$ on $i(K)$ coindides with the push-forward diffeology from $K$ to $i(K).$
%We consider also the corresponding sequence of Lie algebras
%$$ 0 \rightarrow \mathfrak{k} \underrightarrow{i'} \mathfrak{g} \underrightarrow{p}  \mathfrak{h} \rightarrow 0 . $$
%Then, 
%\begin{itemize}
%\item The Lie algebras $\mathfrak{k}$ and $\mathfrak{h}$ are regular if and only if the
%Lie algebra $\mathfrak{g}$ is regular;
%\item The Fr\"olicher Lie groups $K$ and $H$ are regular if and only if the Fr\"olicher Lie group $G$ is regular.
%\end{itemize}

%\end{Theorem}
We mimick and extend the procedure used in \cite{Ma2013}. 
\begin{Theorem} \label{extension}
Let $\mathcal{A} = \bigoplus_{i \in I} \mathcal{A}_i$ be a Fr\"olicher $I-$graded regular algebra. Let $G$ be a regular Fr\"olicher Lie group, acting on $\mathcal{A}$ componentwise. Then,  $G \oplus A$ is a regular Fr\"olicher Lie group.
\end{Theorem}

\noindent
\textbf{Proof.}
Considering the exact sequence 
$$ 0 \rightarrow 1 + \A \rightarrow G \oplus \A \rightarrow G \rightarrow 0$$
there  is a (global) slice $G \rightarrow G \oplus \{0_\A \}$ so that following \cite{Ma2013} Theorem 1.26. \qed

In our work of Lax-type equations, we use the following group from \cite{Ma2013}:Let $M$ be a 
compact manifold without boundary. We denote by  $ \mathcal{F}Cl $ be the space of
 formal classical pseudo-differential operators acting on $C^\infty(M,\R).$ We denote by $\mathcal{F}Cl^*$ the groups of
the units of the algebras $\mathcal{F}Cl$.
Let $q$ be a formal parameter. 
We define the algebra of formal series 
$\mathcal{F}Cl_q = \left\{ \sum_{t \in \N^*} q^k a_k | \forall k \in \N^*, a_k \in \mathcal{F}Cl \right\}.$
This is obviously an algebra, graded by the order (the valuation) into the variable  $q.$ Thus, setting
$ \A_n = \left\{ q^n a_n | a_n \in \mathcal{F}Cl\right\} ,$
we can set $\A = Cl_q(M,E)$ and state the following consequence of Theorem \ref{regulardeformation}:
Let $\mathcal{F}Cl^{0,*}$ be the Lie group of invertible pseudo-differential operators of order 0. This group is known to be a regular Lie group since Omori, but the most efficient proof is actually in \cite{Glo}, to our knowledge.
We remark a short exact sequence of Fr\"olicher Lie groups:
$$ 0 \rightarrow 1 + \mathcal{F}Cl_q \rightarrow \mathcal{F}Cl^{0,*} + \mathcal{F}Cl_q \rightarrow \mathcal{F}Cl^{0,*} \rightarrow 0,$$  
which satisfies the conditions of Theorem \ref{extension}. Thus, we have the following:
\begin{Theorem}The group $1 + \mathcal{F}Cl_q$ is a regular Fr\"olicher Lie group with regular 
Fr\"olicher Lie algebra $\mathcal{F}Cl_q,$ and
$\mathcal{F}Cl^{0,*} + \mathcal{F}Cl_q$ is a regular Fr\"olicher Lie group with Lie algebra $\mathcal{F}Cl^{0} + \mathcal{F}Cl_q.$
\end{Theorem}

 %One could also develop a similar example, 
%which could stand as a generalized version, 
%with log-polyhomogeneous pseudo-differential operators or with other algebras of non classical operators. 
%These examples are not developed here in order to avoid some too long lists of examples 
%constructed in the same spirit.

	\section{On Lax equations and their symmetries} \label{Laxequation}
A PDE is of Lax type if there is a representation of the solutions $u(t,x)\in C^\infty(\mathbb{R} \times M, \mathbb{C})$ in terms of Lax operators, i.e. a smooth map $u \mapsto L(u) \in \mathcal{F}PDO$ (formal pseudo-differential operators) (very often,  $L$ is a differential operator), and another smooth map $u \mapsto P(u) \in \mathcal{F}PDO$ which satisfy a \textbf{Lax equation} such that $u$ is a solution of the initial PDE if and only if the following equation is fulfilled:
\begin{equation}\label{lax1}
\left\{ \begin{array}{lll} \partial_t L(t) & = & \left[ P(t), L(t) \right] \\
 L(0) & = & \hbox{fixed operator (initial value)} \end{array} \right.
\end{equation}
(here and in the sequel, we write $L$ and $P$ instead of $L(u)$ and $P(u)$ when it carries no ambiguity)
The couple (L,P) is called a \textbf{Lax pair.} If the path $P$
is a smooth path of the Lie algebra $\mathfrak{g}$ of a regular Lie group $G,$
if $G$ acts on a Fr\'echet algebra of operators $\mathcal{B}$ that contains $L(0),$ 
the path \begin{equation} \label{solution1} L(t) = Ad_{Exp_GP(t)} L(0) \end{equation}
is a solution of equation \ref{lax1}, yet very often a formal solution. 

\vskip 6pt
\noindent
\textbf{Example: the KdV equation.} The KdV equation reads as $\partial_t u  = 6 u \partial_x u - \partial^3_xu$ where $u(t,x) \in C^\infty(\R^2,\R)$ ad has a Lax pair $L = -\partial_x^2 + u$ and $P = -4\partial^3_x + 3 (\partial_x  u + u \partial_x)$. the operator $P$ is of order 4, so that there is no Lie group $G$  such that $ Exp_GP(t)$ exists. Moreover, P depends on $u$ (essentially because the KdV equation is non linear). 

\vskip 6 pt 
Let us note $\s$ the set of solutions of the initial PDE, which is assumed non empty and 
equipped with the diffeology spanned by (see \cite{Lo1992} for the link between diffeological spaces and Fr\'echet manifolds): 

- the trace diffeology as a subset of  $C^\infty(\mathbb{R} \times M, \mathbb{C}),$ 

- and the pull-back of the diffeology of $\mathcal{F}PDO$ from the maps $u \mapsto L(u)$ and $u \mapsto P(u).$
  
The \textbf{total set of symmetries} of the initial PDE is the group $Diff(\s),$ which is a diffeological group. This space of symmetries is actually, to our knowledge, not studied. Instead of working with $\s,$ we work with $L(\s).$
If $\s$ and $L(\s)$ are (diffeologically) isomorphic, then $Diff(\s)$ and $Diff(L(s))$ are also isomorphic. We now restrict ourselves to smooth linear maps acting on the vector space spanned by $L(\s)$ and $P(\s)$ in $\mathcal{F}PDO. $ Let us write formally the action of such a symmetry S:
if $L(t)$ is a solution of \ref{lax1}, $S(t).L(t)$ is also a solution,
%\begin{eqnarray*} \partial_t(S.L)(t) & = & (\partial_t S).L(t) + S. (\partial_t L) (t)  =  (\partial_t S).L(t) + S.[P,L](t) \\
%& = & (\partial_t S).L(t) + \left( S.[P,L](t) -[P,S.L](t) \right) + \partial_t(S.L)(t)  
%\end{eqnarray*}
from which we get \begin{equation} \label{symmetry1}
(\partial_t S).L(t) = \left[ ad_P,S \right] .L(t)\end{equation}
Here, the map $S$ is a smooth map $S: L(\s) \rightarrow \mathcal{L}(Span(L(\s))).$  The map $\phi : S \mapsto (\partial_t S ).L(t) -  \left[ ad_{P}, S\right] .L(t)$ is linear and the (restricted groups of) symmetries of (\ref{symmetry1}) are the zeros of $\phi.$ 
This relation is now linear in $S$ which allows to pass to infinitesimal symmetries, if the algebra of symmetries under consideration is equipped with the functional diffeology. We get here symmetries that are not in general exactly the ones described in \cite{Olv}, where \textit{projectable} symmetries are the symmetries coming from changing of coordinates, i.e. infinitesimal symmetries in $Vect(\R \times M).$
 
\section{Integration, symmetries and time scaling}
We only assume that both $P$ and $L$ are in a fixed Fr\'echet algebra $\mathcal{A}$ with unit element,
or in a $c^\infty-$algebra if one prefers to work in the convenient setting \cite{KM}. 
Let us now build a corresponding Lax equation in $\mathcal{A}[[q]].$
We consider the paths $P(qt)$ and $L(qt)$ obtained by time scaling $t \mapsto qt. $
Then, $ \partial_t L(qt)  =  q(\partial_tL)(qt) 
=\left[ qP(qt), L(qt) \right] $ 
for a fixed parameter $q.$ We note by $L_q(t)=L(qt)$ and by $P_q(t) = qP(qt).$
We get the following equation:
\begin{equation}\label{lax2}
\left\{ \begin{array}{lll} \partial_t L_q(t) & = & \left[ P_q(t), L_q(t) \right] \\
 L(0) & = & \hbox{fixed operator (initial value) in } \mathcal{A} \end{array} \right.
\end{equation}
Let $val_q$ be the valuation of formal series in $\mathcal{A}[[q]]$ with respect to the $q$ variable.
We remark that $val_qL_q=0$ and $val_qP_q=1. $ We note by $\mathcal{A}[[q]]_{>0}$ the ideal made of formal series S such that $val_qS>0.$
\begin{Theorem}
The solutions of equation \ref{lax2} in $\mathcal{A}[[q]]$ are such that:
$ L_q(t)= exp(P_q)(t).L(0).\left(exp(P_q)(t)\right)^{-1}$
where the the map $exp$ is the group exponential $\mathcal{A}[[q]]_{>0} \rightarrow Id + \mathcal{A}[[q]]_{>0}.$
\end{Theorem}

The proof is a straightforward consequence of basic results on Lie groups. 
The serie $exp(P_q)(t),$ read as 
$  exp(P_q)(t) =  \sum_{i = 0}^\infty a_i(q) $ where 
$a_i (q) = \int_{t\geq s_1 \geq ... \geq s_i\geq 0} \left[ \prod_{j = 1}^i P_q(s_j)\right] (ds)^i
$

Let us now look for symmetries of a Lax equation. A symmetry is a path  $S$ of linear invertible operators on $\mathcal{A}$ such that,  Assuming smoothness, we shall quickly go into more restricted classes of symmetries along the lines of the last section. 
Applying the time scaling, we get, with the obvious notations: 
\begin{equation} \label{symmetry2}
(\partial_t S_q ).L_q(t) =  \left[ ad_{P_q}, S_q\right] .L_q(t)\end{equation}
The map $S \rightarrow S_q$ is an homomorphism from the group of symmetries of (\ref{lax1}) to the group of symmetries of (\ref{lax2}), and it appears to us that there should exist symmetries of (\ref{lax2}) that are not induced from symmetries of (\ref{lax1}).  
The map $\phi_q : S_q \mapsto (\partial_t S_q ).L_q(t) -  \left[ ad_{P_q}, S_q\right] .L_q(t)$ is linear and the symmetries of (\ref{symmetry2}) are the zeros of $\phi_q.$ Such a problem appears non relevant to the methods of resolution of this paper, and we leave the question of solving these two equations open. Let us now turn to a special class of solutions. 

\section{Symmetries as holonomy elements}

Let us now simplify this equation, avoiding the $L_q -$ term. 
Then, we get another Lax-type equation
\begin{equation} \label{symmetry3}
\partial_t S_q  =  \left[ ad_{P_q}, S_q\right] \end{equation}
and we can remark that the operator $ad_{P_q}$ is an inner derivation of $A,$ which is of order 1 in $q$ since $P_q$ is of order 1. Let $In(A)$ be the Lie algebra of inner derivations of $A.$ Let $In_q(A)$ be the $q-$graded algebra of operators spanned by $  qIn_A  ,$ endowed with the push-forward Fr\"olicher structure from $A.$

We have to check:

\begin{enumerate}
\item \label{no1}$In_q$ is a smooth regular algebra
\item \label{no2} $Id_A + In_q(A)$ is a regular Fr\"olicher Lie group with Lie algebra $In_q(A).$ \end{enumerate}

Let us remark that (\ref{no2}) is a straightforward consequence of (\ref{no1}) and Theorem \ref{regulardeformation}.
Now, we recall that smoothness in $In(A)$ is induced by smoothness in $A.$ Moreover, the inclusion $In(A) \rightarrow C^\infty(A,A)$ is smooth in the Fr\"olicher sense \cite{KM}. So that $In_q(A)$ is a smooth algebra,
where the composition is smooth and bilinear. Finally the only checkpoint is that the paths $ad_{a(t)} ad_{b(t)}...$ are integrable . This is in particular true in algebras of formal pseudo-differential operators, using the rules of composition of formal symbols.
We can now apply the procedure that we used for equations (\ref{lax2}): the exponential $exp_{Id_A + In_q(A)} $ exists and
$$ S_q (t) = Exp_{Id_A + In_q(A)}(ad_{P_q}) . S_q(0) .  \left(Exp_{Id_A + In_q(A)}(ad_{P_q})\right)^{-1}$$
is the unique solution to equation (\ref{symmetry2}) with initial value $S_q(0).
$ We now analyze equation \ref{symmetry3}: it extends to the path space of $C^\infty(M,\R)[[q]] \times A_q$ which can be viewed as a trivial vector bundle. Setting $\nabla = d + ad_{P_q}$, we get a smooth connection on this fiber bundle. By the last discussion, $ad_{P_q}$ takes values in a regular Fr\"olicher group, and hence any path on $C^\infty(M)[[q]]$ lifts to a smooth path on $A_q$ by holonomy theorem \cite{Ma2013}, and for any linear map S(0) that transforms an initial value $L(0)$ into another initial solution $S(0).L(0)$, we get a smooth path of operators $t \mapsto S_q(t)$ such that, if $t \mapsto L_q(t)$ is a smooth path which is solution of \ref{lax2}, the path  $t \mapsto S_q(t).L_q(t)$ is also a solution of \ref{lax2}.

\end{document}